\documentclass[12pt]{article}
\usepackage{amsfonts}
\usepackage{amsmath}
\usepackage{graphicx}

\input tcilatex
\input tcilatex
\begin{document}

\author{Aysun Tok Onarcan$^{1}$,Nihat Adar$^{2}$, \.{I}diris Dag$^{2}$ \\
Informatics Department$^{1}$,\\
Computer Engineering Department$^{2}$,\\
Eskisehir Osmangazi University, 26480, Eskisehir, Turkey}
\title{Numerical Solutions of Reaction-Diffusion Equation Systems with
Trigonometric Quintic B-spline Collocation Algorithm}
\maketitle

\begin{abstract}
In this study, the numerical solutions of reaction-diffusion systems are
investigated via the trigonometric quintic B-spline nite element collocation
method. These equations appear in various disciplines in order to describe
certain physical facts, such as pattern formation, autocatalytic chemical
reactions and population dynamics. The Schnakenberg, Gray-Scott and
Brusselator models are special cases of reaction-diffusion systems
considered as numerical examples in this paper. For numerical purposes,
Crank-Nicolson formulae are used for the time discretization and the
resulting system is linearized by Taylor expansion. In the finite element
method, a uniform partition of the solution domain is constructed for the
space discretization. Over the mentioned mesh, dirac-delta function and
trigonometric quintic B-spline functions are chosen as the weighted function
and the bases functions, respectively. Thus, the reaction-diffusion system
turns into an algebraic system which can be represented by a matrix equation
so that the coeffcients are block matrices containing a certain number of
non-zero elements in each row. The method is tested on different problems.
To illustrate the accuracy, error norms are calculated in the linear problem
whereas the relative error is given in other nonlinear problems. Subject to
the character of the nonlinear problems, the occurring spatial patterns are
formed by the trajectories of the dependent variables. The degree of the
base polynomial allows the method to be used in high-order differential
equation solutions. The algorithm produces accurate results even when the
time increment is larger. Therefore, the proposed Trigonometric Quintic
B-spline Collocation method is an effective method which produces acceptable
results for the solutions of reaction-diffusion systems.
\end{abstract}

\section{Introduction}

The reaction diffusion (RD) system is used to model chemical exchange
reactions, the transport of ground water in an aquifer, pattern formation in
the study of biology, chemistry and ecology. The RD system exhibits very
rich dynamics behavior including periodic and quasi-periodic solutions.
Theoretical studies have been developed to describe such dynamic behaviors.
Most reaction-diffusion systems includes the nonlinear reaction term making
it is diffcult to solve analytically. Attempts have been made to look for
the numerical solutions to reveal more dynamic behaviors of the RD system.

The spline functions of various degrees are accompanied to construct
numerical methods to solve di erential equations of certain order, since the
resulting matrix system is always diagonal and can be solved easily and
approximate solutions having the degree accuracy of less than the degree of
the spline functions, can be set up. High order continuous di erentiable
approximate solutions can be produced by way of using high order spline
functions as solutions of the di erential equations. B-splines are de ned as
a basis of the spline space \cite{b21}. Polynomial B-splines are extensively
used for nding numerical solutions of di erential equations, function
approximation and computer-aided design. The numerical procedure based on
the B-spline collocation method has been increasingly applied for nonlinear
evolution equations in various elds of science. However, application of
trigonometric B-spline collocation methods to nonlinear evolution problems
is few in comparison with the collocation method based on polynomial
B-spline functions. The numerical methods for solving types of ordinary di
erential equations with quadratic and cubic trigonometric B-spline are given
by A. Nikolis \cite{b1,b8}. Linear two point boundary value problems of the
order of two are solved using the trigonometric cubic B- spline(TCB)
interpolation method \cite{b16}. Another numerical method employing the TCB
is set up to solve a class of linear two-point singular boundary value
problems in the study \cite{b18}. Recently, a collocation nite di erence
scheme based on the TCB has been developed for the numerical solution of a
one-dimensional hyperbolic equation (wave equation) with a nonlocal
conservation condition \cite{b19}. A new two-time level implicit technique
based on the TCB, is proposed for the approximate solution of a nonclassical
di usion problem with a nonlocal boundary condition in the study \cite{b20}.
A new three-time level implicit approach, based on the TCB is presented for
the approximate solution of the Generalized Nonlinear Klien-Gordon equation
with Dirichlet boundary conditions \cite{kgo}. Some research in the
literature \cite{b15} has established spline-based numerical approaches for
solving reaction-difussion equation systems but without the trigonometric
B-spline, to our knowledge. In this paper, trigonometric quintic
B-splines(TQB) are used to establish a collocation method with suggested
numerical method being applied to nd numerical solutions of a
reaction-diffusion equation system. As a result, the present method makes it
possible to approximate solutions as well as derivatives up to an order of
four at each point of the problem domain.

When reaction-diffusion systems are studied, it can be understood that
different species interact with each other, and also that in chemical
reactions two different chemical substances generate new substances, for
example. For modeling these types of events, which have more than one
dependent variable, differential equation systems have been used.
One-dimensional time-dependent reaction-diffusion equation systems can be
defined as follows:

\begin{equation}
\begin{tabular}{l}
$\dfrac{\partial U}{\partial t}=D_{u}\dfrac{\partial ^{2}U}{\partial x^{2}}%
+F(U,V)$ \\ 
$\dfrac{\partial V}{\partial t}=D_{v}\dfrac{\partial ^{2}V}{\partial x^{2}}%
+G(U,V)$%
\end{tabular}
\label{r3}
\end{equation}%
where $U=U(x,t),V=V(x,t),\Omega \subset R^{2}$ is a problem domain, $D_{u}$
and $D_{v}$ are the diffusion coefficients of $U$ and $V$ respectively, $F$
and $G$ are the growth and interaction functions that represents the
reactions of the system. $F$ and $G$ are always nonlinear functions. A
general one dimensional reaction-diffusion equation system which includes
all models we mentioned in this paper, is expressed as:

\begin{equation}
\begin{tabular}{l}
$\dfrac{\partial U}{\partial t}=a_{1}\dfrac{\partial ^{2}U}{\partial x^{2}}%
+b_{1}U+c_{1}V+d_{1}U^{2}V+e_{1}UV+m_{1}UV^{2}+n_{1}$ \\ 
\\ 
$\dfrac{\partial V}{\partial t}=a_{2}\dfrac{\partial ^{2}V}{\partial x^{2}}%
+b_{2}U+c_{2}V+d_{2}U^{2}V+e_{2}UV+m_{2}UV^{2}+n_{2}$%
\end{tabular}
\label{r1}
\end{equation}%
The solution region of the problem$(-\infty ,\infty )$ should be restricted
as $(x_{0},x_{N})$ for computational purpose. In this case, system (\ref{r1}%
)'s initial conditions are either the homogeny Dirichlet boundary conditions 
\begin{equation}
\begin{tabular}{l}
$U(x_{0},t)=U(x_{N},t)=0,$ \\ 
$V(x_{0},t)=V(x_{N},t)=0,$%
\end{tabular}
\label{r2}
\end{equation}%
or homogeny Neumann boundary conditions%
\begin{equation}
\begin{tabular}{l}
$U_{x}(x_{0},t)=U_{x}(x_{N},t)=0,$ \\ 
$V_{x}(x_{0},t)=V_{x}(x_{N},t)=0$%
\end{tabular}
\label{r4}
\end{equation}%
will be used. Appropriate coefficients of the system (\ref{r1}) for each
test problem will be selected depending on the characteristics of each model
in the following sections and documented in Table 1:

\begin{equation*}
\begin{tabular}{|l|}
\hline
Table 1: The coefficient regulations for model system \\ \hline
\begin{tabular}{lllllllllllllll}
Test Problem & $a_{1}$ & $a_{2}$ & $b_{1}$ & $b_{2}$ & $c_{1}$ & $c_{2}$ & $%
d_{1}$ & $d_{2}$ & $e_{1}$ & $e_{2}$ & $m_{1}$ & $m_{2}$ & $n_{1}$ & $n_{2}$
\\ \hline
Linear & $d$ & $d$ & $-a$ & $0$ & $1$ & $-b$ & $0$ & $0$ & $0$ & $0$ & $0$ & 
$0$ & $0$ & $0$ \\ 
Brusselator & $\varepsilon _{1}$ & $\varepsilon _{2}$ & $-(B+1)$ & $B$ & $0$
& $0$ & $1$ & $-1$ & $0$ & $0$ & $0$ & $0$ & $A$ & $0$ \\ 
Schnakenberg & $1$ & $d$ & $-\gamma $ & $0$ & $0$ & $0$ & $\gamma $ & $%
-\gamma $ & $0$ & $0$ & $0$ & $0$ & $\gamma a$ & $\gamma b$ \\ 
Gray-Scott & $\varepsilon _{1}$ & $\varepsilon _{2}$ & $-f$ & $0$ & $0$ & $%
-(f+k)$ & $0$ & $0$ & $0$ & $0$ & $-1$ & $1$ & $f$ & $0$%
\end{tabular}
\\ \hline
\end{tabular}%
\end{equation*}

\section{The \textbf{Trigonometric Quintic B-spline Collocation Method}}

Consider the solution space of the differential problem $[a=x_{0},b=x_{N}]$
is partitioned into a mesh of uniform length $h=x_{m+1}-x_{m}$ by knots $%
x_{m}$ \ where \ $m=-2,\ldots ,N+2.$ On this partition, together with
additional knots $x_{N-2},x_{N-1},x_{N+1},x_{N+2}$ outside the problem
domain, the trigonometric quintic B-spline $T_{m}^{5}(x)$ \ basis functions
at knots is given by

\begin{equation}
T_{m}^{5}(x)=\frac{1}{\theta }\left \{ 
\begin{tabular}{ll}
$p^{5}(x_{m-3}),$ & $x\in \left[ x_{m-3},x_{m-2}\right] $ \\ \hline
$-p^{4}(x_{m-3})p(x_{m-1})-p^{3}(x_{m-3})p(x_{m})p(x_{m-3})$ &  \\ 
$-p^{2}(x_{m-3})p(x_{m+1})p^{2}(x_{m-2})-p(x_{m-3})p(x_{m+2})p^{3}(x_{m-2})$
&  \\ 
$-p(x_{m+3})p^{4}(x_{m-2}),$ & $x\in \left[ x_{m-2},x_{m-1}\right] $ \\ 
\hline
$p^{3}(x_{m-3})p^{2}(x_{m})+p^{2}(x_{m-3})p(x_{m+1})p(x_{m-2})p(x_{m})$ & 
\\ 
$%
+p^{2}(x_{m-3})p^{2}(x_{m+1})p(x_{m-1})+p(x_{m+3})p(x_{m+2})p^{2}(x_{m-2})p(x_{m}) 
$ &  \\ 
$%
+p(x_{m-3})p(x_{m+2})p(x_{m-2})p(x_{m+1})p(x_{m-1})+p(x_{m-3})p^{2}(x_{m+2})p^{2}(x_{m-1}) 
$ &  \\ 
$%
+p(x_{m+3})p^{3}(x_{m-2})p(x_{m})+p(x_{m+3})p^{2}(x_{m-2})p(x_{m+1})p(x_{m-1}) 
$ &  \\ 
$+p(x_{m+3})p(x_{m-2})p(x_{m+2})p^{2}(x_{m-1})+p^{2}(x_{m+3})p^{3}(x_{m-1}),$
& $x\in \left[ x_{m-1},x_{m}\right] $ \\ \hline
$-p^{2}(x_{m-3})p^{3}(x_{m+1})-p(x_{m-3})p(x_{m+2})p(x_{m-2})p^{2}(x_{m+1})$
&  \\ 
$%
-p(x_{m-3})p^{2}(x_{m+2})p(x_{m-1})p(x_{m+1})-p(x_{m-3})p^{3}(x_{m+2})p(x_{m}) 
$ &  \\ 
$%
-p(x_{m+3})p^{2}(x_{m-2})p^{2}(x_{m})-p(x_{m+3})p(x_{m-2})p(x_{m+2})p(x_{m-1})p(x_{m+1}) 
$ &  \\ 
$-p(x_{m+3})p(x_{m-2})p^{2}(x_{m+2})p(x_{m})-p^{2}(x_{m+3})p^{2}(x_{m-3})$ & 
\\ 
$-p^{2}(x_{m+3})p(x_{m-1})p(x_{m+2})p(x_{m})-p^{3}(x_{m+3})p^{2}(x_{m}),$ & $%
x\in \left[ x_{m},x_{m+1}\right] $ \\ \hline
$%
p(x_{m-3})p^{4}(x_{m+2})+p(x_{m+3})p(x_{m-2})p^{3}(x_{m+2})+p^{2}(x_{m+3})p(x_{m-1})p^{2}(x_{m+2}) 
$ &  \\ 
$+p^{3}(x_{m+3})p(x_{m})p(x_{m+2})+p^{4}(x_{m+3})p(x_{m+1}),$ & $x\in \left[
x_{m+1},x_{m+2}\right] $ \\ \hline
$-p^{5}(x_{m+3}),$ & $x\in \left[ x_{m+2},x_{m+3}\right] $ \\ \hline
$0,$ & $dd$%
\end{tabular}%
\right.  \label{r12}
\end{equation}

where the $p(x_{m}),$ $\Theta $ and$\ m$ are:%
\begin{equation*}
\begin{tabular}{l}
$p(x_{m})=\sin (\frac{x-x_{m}}{2}),$ \\ 
$\Theta =\sin (\frac{5h}{2})\sin (2h)\sin (\frac{3h}{2})\sin (h)\sin (\frac{h%
}{2}),$ \\ 
$m=O(1)N$%
\end{tabular}%
\end{equation*}

The $T_{m}^{5}(x)$ \ functions and its principle derivatives varnish outside
the region $\left[ x_{m-3},x_{m+3}\right] $. The set of those B-splines
\bigskip $T_{m}^{5}(x)$ $,m=-2,...,N+2$ are a basis for the trigonometric
spline space. An approximate solution $U_{N}(x,t)$ and $V_{N}(x,t)$ to the
unknown solution $U(x,t)$ and $V(x,t)$ can be assumed of the forms

\begin{equation}
\begin{tabular}{ll}
$U_{N}(x,t)=\underset{i=-2}{\overset{N+2}{\dsum }}T_{i}^{5}(x)\delta _{i}(t)$
& $V_{N}(x,t)=\underset{i=-2}{\overset{N+2}{\dsum }}T_{i}^{5}(x)\gamma
_{i}(t)$%
\end{tabular}
\label{r14}
\end{equation}%
where $\delta _{i}$ and $\gamma _{i}$ are time dependent parameters to be
determined from the collocation points $x_{i},$ $i=0,...,N$ with boundary
and initial conditions.

$T_{m}^{5}(x)$ trigonometric\ quintic B-spline functions are zero behind the
interval $[x_{m-3},x_{m+3}]$ and $T_{m}^{5}(x)$ functions sequentially
covers six elements in the interval $[x_{m-3},x_{m+3}]$ \ so that, each \ $%
[x_{m},x_{m+1}]$ finite element is covered by the six $%
T_{m-2}^{5},T_{m-1}^{5},T_{m}^{5},T_{m+1}^{5},T_{m+2}^{5},$ and $T_{m+3}^{5}$
trigonometric quintic B-spline. In this case (\ref{r14}) the approach is
given as:

\begin{equation}
\begin{array}{c}
U_{N}(x,t)=\underset{i=m-2}{\overset{m+3}{\dsum }}T_{i}^{5}(x)\delta
_{i}=T_{m-2}^{5}(x)\delta _{m-2}+T_{m-1}^{5}(x)\delta
_{m-1}+T_{m}^{5}(x)\delta _{m}+T_{m+1}^{5}(x)\delta _{m+1} \\ 
\multicolumn{1}{r}{+T_{m+2}^{5}(x)\delta _{m+2}+T_{m+3}^{5}(x)\delta _{m+3}}
\\ 
V_{N}(x,t)=\underset{i=m-2}{\overset{m+3}{\dsum }}T_{i}^{5}(x)\gamma
_{i}=T_{m-2}^{5}(x)\gamma _{m-2}+T_{m-1}^{5}(x)\gamma
_{m-1}+T_{m}^{5}(x)\gamma _{m}+T_{m+1}^{5}(x)\gamma _{m+1} \\ 
\multicolumn{1}{r}{+T_{m+2}^{5}(x)\gamma _{m+2}+T_{m+3}^{5}(x)\gamma _{m+3}}%
\end{array}
\label{r15}
\end{equation}

In these numerical approaches, the approximate solutions at the knots can be
written in terms of the time parametes using $T_{m}^{5}(x)$ and Eq.(\ref{r14}%
). After this, by also making necessary calculations, we can write $%
T_{m}^{5}(x)$ functions for $U_{m}$ and $V_{m}$ and its first, second,third
and fourth derivatives at the knots $x_{m}$ are given in terms of parameters
by the following relationships.

\begin{equation}
\begin{tabular}{l}
\begin{tabular}{l}
$U_{m}=\alpha _{1}\delta _{m-2}+\alpha _{2}\delta _{m-1}+\alpha _{3}\delta
_{m}+\alpha _{2}\delta _{m+1}+\alpha _{1}\delta _{m+2}$ \\ 
$U_{m}^{\prime }=-\alpha _{4}\delta _{m-2}-\alpha _{5}\delta _{m-1}+\alpha
_{5}\delta _{m+1}-\alpha _{4}\delta _{m+2}$ \\ 
$U_{m}^{\prime \prime }=\alpha _{6}\delta _{m-2}+\alpha _{7}\delta
_{m-1}+\alpha _{8}\delta _{m}+\alpha _{7}\delta _{m+1}+\alpha _{6}\delta
_{m+2}$ \\ 
$U_{m}^{\prime \prime \prime }=-\alpha _{9}\delta _{m-2}+\alpha _{10}\delta
_{m-1}-\alpha _{10}\delta _{m+1}-\alpha _{9}\delta _{m+2}$ \\ 
$U_{m}^{\prime \prime \prime \prime }=\alpha _{11}\delta _{m-2}+\alpha
_{12}\delta _{m-1}+\alpha _{13}\delta _{m}+\alpha _{12}\delta _{m+1}+\alpha
_{11}\delta _{m+2}$%
\end{tabular}
\\ 
\begin{tabular}{l}
$V_{m}=\alpha _{1}\gamma _{m-2}+\alpha _{2}\gamma _{m-1}+\alpha _{3}\gamma
_{m}+\alpha _{2}\gamma _{m+1}+\alpha _{1}\gamma _{m+2}$ \\ 
$V_{m}^{\prime }=-\alpha _{4}\gamma _{m-2}-\alpha _{5}\gamma _{m-1}+\alpha
_{5}\gamma _{m+1}+\alpha _{4}\gamma _{m+2}$ \\ 
$V_{m}^{\prime \prime }=\alpha _{6}\gamma _{m-2}+\alpha _{7}\gamma
_{m-1}+\alpha _{8}\gamma _{m}+\alpha _{7}\gamma _{m+1}+\alpha _{6}\gamma
_{m+2}$ \\ 
$V_{m}^{\prime \prime \prime }=-\alpha _{9}\gamma _{m-2}+\alpha _{10}\gamma
_{m-1}-\alpha _{10}\gamma _{m+1}+\alpha _{9}\gamma _{m+2}$ \\ 
$V_{m}^{\prime \prime \prime \prime }=\alpha _{11}\gamma _{m-2}+\alpha
_{12}\gamma _{m-1}+\alpha _{13}\gamma _{m}+\alpha _{12}\gamma _{m+1}+\alpha
_{11}\gamma _{m+2}$%
\end{tabular}%
\end{tabular}
\label{r16}
\end{equation}%
where the coefficients are:

\begin{equation*}
\begin{array}{l}
\alpha _{1}=\dfrac{\sin ^{5}(\frac{h}{2})}{\Theta } \\ 
\\ 
\alpha _{2}=\dfrac{2\sin ^{5}(\frac{h}{2})\cos (\frac{h}{2})(16\cos ^{2}(%
\frac{h}{2})-3)}{\Theta } \\ 
\\ 
\alpha _{3}=\dfrac{2(1+48\cos {}^{4}(\frac{h}{2})-16\cos ^{2}(\frac{h}{2}%
)\sin {}^{5}(\frac{h}{2}))}{\Theta } \\ 
\\ 
\alpha _{4}=\dfrac{\frac{5}{2}\sin ^{4}(\frac{h}{2})\cos (\frac{h}{2})}{%
\Theta } \\ 
\\ 
\alpha _{5}=\dfrac{5\sin ^{4}(\frac{h}{2})\cos ^{2}(\frac{h}{2})(8\cos ^{2}(%
\frac{h}{2})-3)}{\Theta } \\ 
\\ 
\alpha _{6}=\dfrac{\frac{5}{4}\sin ^{3}(\frac{h}{2})(5\cos ^{2}(\frac{h}{2}%
)-1)}{\Theta } \\ 
\\ 
\alpha _{7}=\dfrac{\frac{5}{2}\sin ^{3}(\frac{h}{2})(\cos (\frac{h}{2}%
)(-15\cos ^{2}(\frac{h}{2})+3+16\cos ^{4}(\frac{h}{2}))}{\Theta } \\ 
\\ 
\alpha _{8}=\dfrac{-\frac{5}{2}\sin ^{3}(\frac{h}{2})(16\cos ^{6}(\frac{h}{2}%
)-5\cos ^{6}(\frac{h}{2})+1)}{\Theta } \\ 
\\ 
\alpha _{9}=\dfrac{\frac{5}{8}\sin ^{2}(\frac{h}{2})\cos (\frac{h}{2}%
)(25\cos ^{2}(\frac{h}{2})-13)}{\Theta } \\ 
\\ 
\alpha _{10}=\dfrac{-\frac{5}{4}\sin ^{2}(\frac{h}{2})(\cos ^{2}(\frac{h}{2}%
)(8\cos ^{4}(\frac{h}{2})-35\cos ^{2}(\frac{h}{2})+15)}{\Theta } \\ 
\alpha _{11}=\dfrac{\frac{5}{16}(125\cos ^{4}(\frac{h}{2})-114\cos ^{2}(%
\frac{h}{2})+13)\sin (\frac{h}{2}))}{\Theta } \\ 
\\ 
\alpha _{12}=\dfrac{-\frac{5}{8}\sin (\frac{h}{2})\cos (\frac{h}{2})(176\cos
^{6}(\frac{h}{2})-137\cos ^{4}(\frac{h}{2})-6\cos ^{2}(\frac{h}{2})+15)}{%
\Theta } \\ 
\\ 
\alpha _{13}=\dfrac{\frac{5}{8}(92\cos ^{6}(\frac{h}{2})-117\cos ^{4}(\frac{h%
}{2})+62\cos ^{2}(\frac{h}{2})-13)(-1+4\cos ^{2}(\frac{h}{2})\sin (\frac{h}{2%
}))}{\Theta }%
\end{array}%
\end{equation*}

The Crank--Nicholson scheme 
\begin{equation}
\begin{tabular}{ll}
$U_{t}=\dfrac{U^{n+1}-U^{n}}{\Delta t},$ & $U=\dfrac{U^{n+1}+U^{n}}{2}$ \\ 
$V_{t}=\dfrac{V^{n+1}-V^{n}}{\Delta t},$ & $V=\dfrac{V^{n+1}+V^{n}}{2}$%
\end{tabular}
\label{r17}
\end{equation}%
is used to discretize time variables of the unknown $U$ and $V$ and their
derivatives, to have the time integrated reaction-difussion equation system:

\begin{equation}
\begin{tabular}{r}
$\dfrac{U^{n+1}-U^{n}}{\Delta t}-a_{1}\dfrac{U_{xx}^{n+1}+U_{xx}^{n}}{2}%
-b_{1}\dfrac{U^{n+1}+U^{n}}{2}-c_{1}\dfrac{V^{n+1}+V^{n}}{2}-d_{1}\dfrac{%
(U^{2}V)^{n+1}+(U^{2}V)^{n}}{2}$ \\ 
$-e_{1}\dfrac{(UV)^{n+1}+(UV)^{n}}{2}-m_{1}\dfrac{(UV^{2})^{n+1}+(UV^{2})^{n}%
}{2}-n_{1}=0$ \\ 
$\dfrac{V^{n+1}-V^{n}}{\Delta t}-a_{2}\dfrac{V_{xx}^{n+1}+V_{xx}^{n}}{2}%
-b_{2}\dfrac{U^{n+1}+U^{n}}{2}-c_{2}\dfrac{V^{n+1}+V^{n}}{2}-d_{2}\dfrac{%
(U^{2}V)^{n+1}+(U^{2}V)^{n}}{2}$ \\ 
$-e_{2}\dfrac{(UV)^{n+1}+(UV)^{n}}{2}-m_{2}\dfrac{(UV^{2})^{n+1}+(UV^{2})^{n}%
}{2}-n_{2}=0$%
\end{tabular}
\label{r18}
\end{equation}%
where $U^{n+1}=U(x,t)$ \ and\ $V^{n+1}=V(x,t)$ are the solutions of the
equations at the $(n+1)$th time level. Here $t_{n+1}=t_{n}+\Delta t$ and $%
\Delta t$ is the time step, superscripts denote the $n$ th level $%
t_{n}=n\Delta t.$

The nonlinear terms $(U^{2}V)^{n+1},$ $(UV^{2})^{n+1}$and $(UV)^{n+1}$ in
equation (\ref{r18}) is linearized by using the following forms \ref{r19} .

\begin{equation}
\begin{tabular}{l}
$%
(U^{2}V)^{n+1}=U^{n+1}U^{n}V^{n}+U^{n}U^{n+1}V^{n}+U^{n}U^{n}V^{n+1}-2U^{n}U^{n}V^{u} 
$ \\ 
$%
(UV^{2})^{n+1}=U^{n+1}V^{n}V^{n}+U^{n}V^{n+1}V^{n}+U^{n}V^{n}V^{n+1}-2U^{n}V^{n}V^{u} 
$ \\ 
$(UV)^{n+1}=U^{n+1}V^{n}+U^{n}V^{n+1}-U^{n}V^{n}$%
\end{tabular}
\label{r19}
\end{equation}%
When we substitute (\ref{r19}) in (\ref{r18}), the linearized general model
equation system takes the form as shown below,

\begin{eqnarray}
-\dfrac{a_{1}}{2}U_{xx}^{n+1}+\beta _{m1}U^{n+1}+\beta _{m2}V^{n+1} &=&%
\dfrac{a_{1}}{2}U_{xx}^{n}+\beta _{m3}U^{n}+\beta _{m4}V^{n}+n_{1}
\label{r20} \\
-\dfrac{a_{2}}{2}V_{xx}^{n+1}+\beta _{m5}U^{n+1}+\beta _{m6}V^{n+1} &=&%
\dfrac{a_{2}}{2}V_{xx}^{n}+\beta _{m7}U^{n}+\beta _{m8}V^{n}+n_{2}  \notag
\end{eqnarray}%
where%
\begin{equation*}
\begin{tabular}{l}
$\beta _{m1}=\dfrac{1}{\Delta t}-\dfrac{b_{1}}{2}-d_{1}U^{n}V^{n}-\dfrac{%
e_{1}}{2}Vn-\dfrac{m_{1}}{2}(V^{n})^{2}$ \\ 
$\beta _{m2}=\dfrac{1}{\Delta t}-\dfrac{c_{1}}{2}-\dfrac{d_{1}}{2}%
(U^{n})^{2}-\dfrac{e_{1}}{2}Un-m_{1}U^{n}V^{n}$ \\ 
$\beta _{m3}=\dfrac{1}{\Delta t}+\dfrac{b1}{2}-\dfrac{m1}{2}(V^{n})^{2}$ \\ 
$\beta _{m4}=\dfrac{c1}{2}-\dfrac{d1}{2}(U^{n})^{2}$ \\ 
$\beta _{m5}=-\dfrac{b_{2}}{2}-d_{2}U^{n}V^{n}-\dfrac{e_{2}}{2}Vn-\dfrac{%
m_{2}}{2}(V^{n})^{2}$ \\ 
$\beta _{m6}=\dfrac{1}{\Delta t}-\dfrac{c_{2}}{2}-\dfrac{d_{2}}{2}%
(U^{n})^{2}-\dfrac{e_{2}}{2}Un-m_{2}U^{n}V^{n}$ \\ 
$\beta _{m7}=\dfrac{b_{2}}{2}-\dfrac{m_{2}}{2}(V^{n})^{2}$ \\ 
$\beta _{m8}=\dfrac{1}{\Delta t}+\dfrac{c_{2}}{2}-\dfrac{d_{2}}{2}%
(U^{n})^{2}.$%
\end{tabular}%
\end{equation*}

To discrete the model system (\ref{r1}) fully by space respectively, we
substitute the approximate solution (\ref{r16}) into (\ref{r20}) yielding
the fully-discretized equations.

\begin{eqnarray}
&&%
\begin{tabular}{l}
$\nu _{m1}\delta _{m-2}^{n+1}+\nu _{m2}\gamma _{m-2}^{n+1}+\nu _{m3}\delta
_{m-1}^{n+1}+\nu _{m4}\gamma _{m-1+}^{n+1}+\nu _{m5}\delta _{m}^{n+1}+\nu
_{m6}\gamma _{m}^{n+1}+$ \\ 
$\nu _{m7}\delta _{m+1}^{n+1}+\nu _{m8}\gamma _{m+1}^{n+1}+\nu _{m9}\delta
_{m+2}^{n+1}+\nu _{m10}\gamma _{m+2}^{n+1}=$ \\ 
$\nu _{m11}\delta _{m-2}^{n}+\nu _{m12}\gamma _{m-2}^{n}+\nu _{m13}\delta
_{m-1}^{n}+\nu _{m14}\gamma _{m-1}^{n}+\nu _{m15}\delta _{m}^{n}+\nu
_{m16}\gamma _{m}^{n}+$ \\ 
$\  \  \  \  \  \  \  \  \  \  \  \  \  \  \  \  \  \  \  \  \  \  \  \  \  \  \nu _{m17}\delta
_{m+1}^{n}+\nu _{m18}\gamma _{m+1}^{n}+\nu _{m19}\delta _{m+2}^{n}+\nu
_{m20}\gamma _{m+2}^{n}+n_{1}$%
\end{tabular}
\label{r21} \\
&&%
\begin{tabular}{l}
$\nu _{m21}\delta _{m-2}^{n+1}+\nu _{m22}\gamma _{m-2}^{n+1}+\nu
_{m23}\delta _{m-1}^{n+1}+\nu _{m24}\gamma _{m-1}^{n+1}+\nu _{m25}\delta
_{m}^{n+1}+\nu _{m26}\gamma _{m}^{n+1}+$ \\ 
$\nu _{m27}\delta _{m+1}^{n+1}+\nu _{m28}\gamma _{m+1}^{n+1}+\nu
_{m29}\delta _{m+2}^{n+1}+\nu _{m30}\gamma _{m+2}^{n+1}=$ \\ 
$\nu _{m31}\delta _{m-2}^{n}+\nu _{m32}\gamma _{m-2}^{n}+\nu _{m33}\delta
_{m-1}^{n}+\nu _{m34}\gamma _{m-1}^{n}+\nu _{m35}\delta _{m}^{n}+\nu
_{m36}\gamma _{m}^{n}+$ \\ 
$\  \  \  \  \  \  \  \  \  \  \  \  \  \  \  \  \  \  \  \  \  \  \  \  \  \  \  \  \  \  \nu
_{m37}\delta _{m+1}^{n}+\nu _{m38}\gamma _{m+1}^{n}+\nu _{m39}\delta
_{m+2}^{n}+\nu _{m40}\gamma _{m+2}^{n}+n_{2}$%
\end{tabular}
\notag
\end{eqnarray}%
where the $\nu _{m}$ coefficients are:

\begin{equation}
\begin{tabular}{llll}
$\nu _{m1}=\beta _{m1}\alpha _{1}-\dfrac{a_{1}}{2}\alpha _{6}$ & $\nu
_{m21}=\beta _{m5}\alpha _{1}$ & $\nu _{m11}=\beta _{m3}\alpha _{1}+\dfrac{%
a_{1}}{2}\alpha _{6}$ & $\nu _{m31}=\beta _{m7}\alpha _{1}$ \\ 
$\nu _{m2}=\beta _{m2}\alpha _{1}$ & $\nu _{m22}=\beta _{m6}\alpha _{1}+%
\dfrac{a_{2}}{2}\alpha _{6}$ & $\nu _{m12}=\beta _{m4}\alpha _{1}$ & $\nu
_{m32}=\beta _{m8}\alpha _{1}-\dfrac{a_{2}}{2}\alpha _{6}$ \\ 
$\nu _{m3}=\beta _{m1}\alpha _{2}-\dfrac{a_{1}}{2}\alpha _{7}$ & $\nu
_{m23}=\beta _{m5}\alpha _{2}$ & $\nu _{m13}=\beta _{m3}\alpha _{2}+\dfrac{%
a_{1}}{2}\alpha _{7}$ & $\nu _{m33}=\beta _{m7}\alpha _{2}$ \\ 
$\nu _{m4}=\beta _{m2}\alpha _{2}$ & $\nu _{m24}=\beta _{m6}\alpha _{2}+%
\dfrac{a_{2}}{2}\alpha _{7}$ & $\nu _{m14}=\beta _{m4}\alpha _{2}$ & $\nu
_{m34}=\beta _{m8}\alpha _{2}-\dfrac{a_{2}}{2}\alpha _{7}$ \\ 
$\nu _{m5}=\beta _{m1}\alpha _{3}-\dfrac{a_{1}}{2}\alpha _{8}$ & $\nu
_{m25}=\beta _{m5}\alpha _{3}$ & $\nu _{m15}=\beta _{m3}\alpha _{3}+\dfrac{%
a_{1}}{2}\alpha _{8}$ & $\nu _{m35}=\beta _{m7}\alpha _{3}$ \\ 
$\nu _{m6}=\beta _{m2}\alpha _{3}$ & $\nu _{m26}=\beta _{m6}\alpha _{3}+%
\dfrac{a_{2}}{2}\alpha _{8}$ & $\nu _{m16}=\beta _{m4}\alpha _{3}$ & $\nu
_{m36}=\beta _{m8}\alpha _{3}-\dfrac{a_{2}}{2}\alpha _{8}$ \\ 
$\nu _{m7}=\beta _{m1}\alpha _{2}-\dfrac{a_{1}}{2}\alpha _{7}$ & $\nu
_{m27}=\beta _{m5}\alpha _{2}$ & $\nu _{m17}=\beta _{m3}\alpha _{2}+\dfrac{%
a_{1}}{2}\alpha _{7}$ & $\nu _{m37}=\beta _{m7}\alpha _{2}$ \\ 
$\nu _{m8}=\beta _{m2}\alpha _{2}$ & $\nu _{m28}=\beta _{m6}\alpha _{2}+%
\dfrac{a_{2}}{2}\alpha _{7}$ & $\nu _{m18}=\beta _{m4}\alpha _{2}$ & $\nu
_{m38}=\beta _{m8}\alpha _{2}-\dfrac{a_{2}}{2}\alpha _{7}$ \\ 
$\nu _{m9}=\beta _{m1}\alpha _{1}-\dfrac{a_{1}}{2}\alpha _{6}$ & $\nu
_{m29}=\beta _{m5}\alpha _{1}$ & $\nu _{m19}=\beta _{m3}\alpha _{1}+\dfrac{%
a_{1}}{2}\alpha _{6}$ & $\nu _{m39}=\beta _{m7}\alpha _{1}$ \\ 
$\nu _{m10}=\beta _{m2}\alpha _{1}$ & $\nu _{m30}=\beta _{m6}\alpha _{1}+%
\dfrac{a_{2}}{2}\alpha _{6}$ & $\nu _{m20}=\beta _{m4}\alpha _{1}$ & $\nu
_{m40}=\beta _{m8}\alpha _{1}-\dfrac{a_{2}}{2}\alpha _{6}$%
\end{tabular}
\label{r22}
\end{equation}

The system (\ref{r21}) can be converted into the following matrix system:

\begin{equation}
A\mathbf{x}^{n+1}=B\mathbf{x}^{n}+F  \label{r23}
\end{equation}

\begin{equation}
\begin{tabular}{l}
$A=\left[ 
\begin{array}{cccccccccccccc}
\nu _{m1} & \nu _{m2} & \nu _{m3} & \nu _{m4} & \nu _{m5} & \nu _{m6} & \nu
_{m7} & \nu _{m8} & \nu _{m9} & \nu _{m10} &  &  &  &  \\ 
\nu _{m21} & \nu _{m22} & \nu _{m23} & \nu _{m24} & \nu _{m25} & \nu _{m26}
& \nu _{m27} & \nu _{m28} & \nu _{m29} & \nu _{m30} &  &  &  &  \\ 
&  & \nu _{m1} & \nu _{m2} & \nu _{m3} & \nu _{m4} & \nu _{m5} & \nu _{m6} & 
\nu _{m7} & \nu _{m8} & \nu _{m9} & \nu _{m10} &  &  \\ 
&  & \nu _{m21} & \nu _{m22} & \nu _{m23} & \nu _{m24} & \nu _{m25} & \nu
_{m26} & \nu _{m27} & \nu _{m28} & \nu _{m29} & \nu _{m30} &  &  \\ 
&  &  & ... & ... & ... & ... & ... & ... & ... & ... & ... & ... & ... \\ 
&  &  &  & \nu _{m1} & \nu _{m2} & \nu _{m3} & \nu _{m4} & \nu _{m5} & \nu
_{m6} & \nu _{m7} & \nu _{m8} & \nu _{m9} & \nu _{m10} \\ 
&  &  &  & \nu _{m21} & \nu _{m22} & \nu _{m23} & \nu _{m24} & \nu _{m25} & 
\nu _{m26} & \nu _{m27} & \nu _{m28} & \nu _{m29} & \nu _{m30}%
\end{array}%
\right] $ \\ 
\\ 
$B=\left[ 
\begin{array}{cccccccccccccc}
\nu _{m11} & \nu _{m12} & \nu _{m13} & \nu _{m14} & \nu _{m15} & \nu _{m16}
& \nu _{m17} & \nu _{m18} & \nu _{m19} & \nu _{m20} &  &  &  &  \\ 
\nu _{m31} & \nu _{m32} & \nu _{m33} & \nu _{m34} & \nu _{m35} & \nu _{m36}
& \nu _{m37} & \nu _{m38} & \nu _{m39} & \nu _{m40} &  &  &  &  \\ 
&  & \nu _{m11} & \nu _{m12} & \nu _{m13} & \nu _{m14} & \nu _{m15} & \nu
_{m16} & \nu _{m17} & \nu _{m18} & \nu _{m19} & \nu _{m20} &  &  \\ 
&  & \nu _{m31} & \nu _{m32} & \nu _{m33} & \nu _{m34} & \nu _{m35} & \nu
_{m36} & \nu _{m37} & \nu _{m38} & \nu _{m39} & \nu _{m40} &  &  \\ 
&  &  & ... & ... & ... & ... & ... & ... & ... & ... & ... & ... & ... \\ 
&  &  &  & \nu _{m11} & \nu _{m12} & \nu _{m13} & \nu _{m14} & \nu _{m15} & 
\nu _{m16} & \nu _{m17} & \nu _{m18} & \nu _{m19} & \nu _{m20} \\ 
&  &  &  & \nu _{m31} & \nu _{m32} & \nu _{m33} & \nu _{m34} & \nu _{m35} & 
\nu _{m36} & \nu _{m37} & \nu _{m38} & \nu _{m39} & \nu _{m40}%
\end{array}%
\right] $%
\end{tabular}
\label{r24}
\end{equation}%
The system (\ref{r24}) is consists of a $2N+2$ linear equation in $2N+10$
unknown parameters with $x^{n+1},x^{n}$ and $F$ being the vectors as shown
below:

\begin{eqnarray*}
\mathbf{x}^{n+1} &=&[\delta _{-2}^{n+1},\gamma _{-2}^{n+1},\delta
_{-1}^{n+1},\gamma _{-1}^{n+1},\delta _{0}^{n+1},\gamma _{0}^{n+1}...,\delta
_{N+1}^{n+1},\gamma _{N+1}^{n+1},\delta _{N+2}^{n+1},\gamma _{N+2}^{n+1}]^{T}
\\
\mathbf{x}^{n} &=&[\delta _{-2}^{n},\gamma _{-2}^{n},\delta _{-1}^{n},\gamma
_{-1}^{n},\delta _{0}^{n},\gamma _{0}^{n}...,\delta _{N+1}^{n},\gamma
_{N+1}^{n},\delta _{N+2}^{n},\gamma _{N+2}^{n}]^{T} \\
F &=&[n_{1},n_{2},n_{1},n_{2},,,n_{1},n_{2}]^{T}
\end{eqnarray*}

To obtain a unique solution an additional eight constraints are needed.\
While $m=0$ \ and $m=N$ by imposing the Dirichlet boundary conditions or the
Neumann boundary conditions this will lead us to new relationships to
eliminate parameters

$\delta _{-2},$ $\delta _{-1},\delta _{N+1},\delta _{N+2},\gamma
_{-2},\gamma _{-1},\gamma _{N+1},\gamma _{N+2}$ from the system (\ref{r23}).
When we eliminate these parameters the resulting $(2N+2)\times (2N+2)$
matrix system can be solved by the Gauss elimination algorithm.

The initial parameters of $\mathbf{x}^{0}=(\delta _{-2}^{0},\gamma
_{-2}^{0},\delta _{-1}^{0},\gamma _{-1}^{0},\delta _{0}^{0},\gamma
_{0}^{0}...,\delta _{N+1}^{0},\gamma _{N+1}^{0},\delta _{N+2}^{0},\gamma
_{N+2}^{0})$ \ must be found to start the iteration process by using both
initial and boundary conditions. The recurrence relationship (\ref{r23})
gives the time evolution of vector $\mathbf{x}^{n}$. Thus, the nodal values $%
U_{N}(x,t)$ and $V_{N}(x,t)$ can be computed via the equations (\ref{r16})
at the knots.

\subsection{Results of The Numerical Solutions}

In this section, we will compare the efficiency and accuracy of the
suggested method on the given reaction-diffusion equation system models. The
obtained results for each model will compare with \cite{b15} and \cite{b3}.
The accuracy of the schemes is measured in terms of the following discrete
error norm

$L_{2}$ $=|U-U_{N}|_{2}$=$\sqrt{h\sum_{j=0}^{N}(U_{j}-(U_{N})_{j}^{n})}$and $%
L_{\infty }=|U-U_{N}|_{\infty }=\underset{j}{\max }|U_{j}-(U_{N})_{j}^{n}|$.

The relative error $=\sqrt{\dfrac{\sum_{j=0}^{N}|U_{j}^{n+1}-U_{j}^{n}|^{2}}{%
\sum_{j=0}^{N}|U_{j}^{n+1}|}}$ is used to measure errors of solutions of the
reaction-diffusion systems that do not have an analytic solution.

\subsubsection{Linear Problem}

It is stated that the terms $F(U,V)$ and $G(U,V)$ are always nonlinear in
the system (\ref{r3}). However, it is not possible to calculate error norms
because of the limitations of the analytical solutions of the nonlinear
system. The linear problem has been solved to examine error norms for
testing this method:

\begin{equation}
\begin{tabular}{l}
$\dfrac{\partial U}{\partial t}=d\dfrac{\partial ^{2}U}{\partial x^{2}}-aU+V$
\\ 
$\dfrac{\partial V}{\partial t}=d\dfrac{\partial ^{2}V}{\partial x^{2}}-bV.$%
\end{tabular}
\label{r6}
\end{equation}

The given equation system described above is a linear reaction-diffusion
system, which has analytical solutions given as:

\begin{equation}
\begin{tabular}{l}
$U(x,t)=(e^{-(a+d)t}+e^{-(b+d)t})\cos (x),$ \\ 
$V(x,t)=(a-b)(e^{-(b+d)t})\cos (x).$%
\end{tabular}
\label{r7}
\end{equation}

Solutions were obtained by solving the reaction-diffusion system (\ref{r6})
in this section. Three different cases were considered in numerical
computation of coefficients in thesystem (\ref{r6}). This system's initial
conditions can be obtained, when $t=0$ in (\ref{r7}) the solutions. When a
solution region is selected as $(0,\dfrac{\pi }{2})$ interval, the boundary
conditions are described as:

\begin{equation}
\begin{tabular}{ll}
$U_{x}(0,t)=0$ & $U(\pi /2,t)=0,$ \\ 
$V_{x}(0,t)=0$ & $V(\pi /2,t)=0.$%
\end{tabular}
\label{r27}
\end{equation}

In numerical calculations, the programme is going to run up to time $t=1$
for various $N$ and $\Delta t$ and the reaction and diffusion mechanism is
examined for different selections of constants $a,b,$and $d.$ The error
values $L_{2}$ and $L_{\infty }$ that have emerged in the solution, are
presented in the tables.

Firstly, the equation system (\ref{r6}) coefficients are chosen as $a=0.1,$ $%
b=0.01$ and $d=1$ which is a diffusion dominated case. The boundary and
initial conditions are chosen to coincide with the polynomial quintic
B-spline collocation method (PQBCM) \cite{b15}. The programme is run up to $%
t=1$ and the obtained results for $U,$ in terms of $L_{2}$ and $L_{\infty }$
norms are given in Table 3.

In Table 3, $L_{2}$ and $L_{\infty }$error norms are calculated for both $U$
and $V,$ for $N=512$ and various $\Delta t$ with results of \cite{b15} and 
\cite{b3} is also given in the same table. When Table 3 is examined, it
seems that, the accuracy of the obtained results for function $V$ are more
efficient than obtained results for function $U$. When we compare the
results, the proposed method has better accuracy aganist the other
references under the same conditions.

\begin{equation*}
\begin{tabular}{|l|l|l|l|l|l|l|l|l|l|}
\hline
\multicolumn{10}{|l|}{Table 3: Error norms $L_{2}$ and $L_{\infty }$ for
diffusion dominant case for $a=0.1,$ $b=0.01,$ $d=1$} \\ \hline
\multicolumn{6}{|l}{TQB} & \multicolumn{4}{|l|}{Polynomial quintic B-spline,%
\cite{b15}} \\ \hline
$N$ & $\Delta t$ & $U$ &  & $V$ &  & $U$ &  & $V$ &  \\ \hline
&  & $L_{2}\times 10^{4}$ & $L_{\infty }\times 10^{4}$ & $L_{2}\times 10^{6}$
& $L_{\infty }\times 10^{6}$ & $L_{2}\times 10^{4}$ & $L_{\infty }\times
10^{4}$ & $L_{2}\times 10^{6}$ & $L_{\infty }\times 10^{6}$ \\ \hline
$512$ & $0.005$ & \multicolumn{1}{|r|}{$0,008090$} & \multicolumn{1}{|r|}{$%
0,009120$} & \multicolumn{1}{|r|}{$0,029344$} & \multicolumn{1}{|r|}{$%
0,033079$} & $0,015123$ & $0,017048$ & $0,062416$ & $0,070361$ \\ \hline
& $0.01$ & \multicolumn{1}{|r|}{$0,053460$} & \multicolumn{1}{|r|}{$0,060265$%
} & \multicolumn{1}{|r|}{$0,216594$} & \multicolumn{1}{|r|}{$0,244162$} & $%
0,060493$ & $0,068193$ & $0,249667$ & $0,281444$ \\ \hline
& $0.02$ & \multicolumn{1}{|r|}{$0,234949$} & \multicolumn{1}{|r|}{$0,264853$%
} & \multicolumn{1}{|r|}{$0,965627$} & \multicolumn{1}{|r|}{$1,088530$} & $%
0,241983$ & $0,272782$ & $0,998702$ & $1,125815$ \\ \hline
& $0.04$ & \multicolumn{1}{|r|}{$0,961033$} & \multicolumn{1}{|r|}{$1,083353$%
} & \multicolumn{1}{|r|}{$3,962253$} & \multicolumn{1}{|r|}{$4,466566$} & $%
0,968068$ & $1,091283$ & $3,995334$ & $4,503855$ \\ \hline
\multicolumn{10}{|l|}{CN-MG method \cite{b3}} \\ \hline
$512$ & $0.005$ &  & $0.0116$ &  &  &  &  &  &  \\ \hline
& $0.01$ &  & $0.0627$ &  &  &  &  &  &  \\ \hline
& $0.02$ &  & $0.267$ &  &  &  &  &  &  \\ \hline
& $0.04$ &  & $1.09$ &  &  &  &  &  &  \\ \hline
\end{tabular}%
\end{equation*}

Secondly, the constants of system equation (\ref{r6}) are selected as $%
a=2,b=1,d=0.001$ which is a reaction dominated case. The programme is run up
to $t=1,$ and the obtained results in terms of $L_{2}$ and $L_{\infty }$
norms are given in Table 4.

In Table 4, $L_{2}$ and $L_{\infty }$ error norms are calculated both for $U$
and $V,$ for $N=512$ and various $\Delta t$ and the results of \cite{b15}
and \cite{b3} are given in the same table.

\begin{equation*}
\begin{tabular}{|l|l|l|l|l|l|l|l|l|l|}
\hline
\multicolumn{10}{|l|}{Table 4: Error norms $L_{2}$ and $L_{\infty }$ for
reaction dominated case for $a=2,$ $b=1,d=0.001$} \\ \hline
\multicolumn{6}{|l}{TQB} & \multicolumn{4}{|l|}{Polynomial quintic B-spline,%
\cite{b15}} \\ \hline
$N$ & $\Delta t$ & $U$ &  & $V$ &  & $U$ &  & $V$ &  \\ \hline
&  & $L_{2}\times 10^{4}$ & $L_{\infty }\times 10^{4}$ & $L_{2}\times 10^{5}$
& $L_{\infty }\times 10^{5}$ & $L_{2}\times 10^{4}$ & $L_{\infty }\times
10^{4}$ & $L_{2}\times 10^{3}$ & $L_{\infty }\times 10^{3}$ \\ \hline
$512$ & $0.005$ & \multicolumn{1}{|r|}{$0,026827$} & \multicolumn{1}{|r|}{$%
0,030241$} & \multicolumn{1}{|r|}{$0,068087$} & \multicolumn{1}{|r|}{$%
0,076753$} & \multicolumn{1}{|r|}{$0,026832$} & \multicolumn{1}{|r|}{$%
0,030247$} & \multicolumn{1}{|r|}{$0,068124$} & \multicolumn{1}{|r|}{$%
0,076795$} \\ \hline
& $0.01$ & \multicolumn{1}{|r|}{$0,107324$} & \multicolumn{1}{|r|}{$0,120984$%
} & \multicolumn{1}{|r|}{$0,272462$} & \multicolumn{1}{|r|}{$0,307141$} & 
\multicolumn{1}{|r|}{$0,107329$} & \multicolumn{1}{|r|}{$0,120989$} & 
\multicolumn{1}{|r|}{$0,272499$} & \multicolumn{1}{|r|}{$0,307183$} \\ \hline
& $0.02$ & \multicolumn{1}{|r|}{$0,429339$} & \multicolumn{1}{|r|}{$0,483984$%
} & \multicolumn{1}{|r|}{$1,089996$} & \multicolumn{1}{|r|}{$1,228729$} & 
\multicolumn{1}{|r|}{$0,429344$} & \multicolumn{1}{|r|}{$0,483990$} & 
\multicolumn{1}{|r|}{$1,090033$} & \multicolumn{1}{|r|}{$1,228771$} \\ \hline
& $0.04$ & \multicolumn{1}{|r|}{$1,717837$} & \multicolumn{1}{|r|}{$1,936481$%
} & \multicolumn{1}{|r|}{$4,360663$} & \multicolumn{1}{|r|}{$4,915683$} & 
\multicolumn{1}{|r|}{$1,717842$} & \multicolumn{1}{|r|}{$1,936487$} & 
\multicolumn{1}{|r|}{$4,360700$} & \multicolumn{1}{|r|}{$4,915725$} \\ \hline
\multicolumn{10}{|l|}{CN-MG method \cite{b3}} \\ \hline
$512$ & $0.005$ &  & $0.0302$ &  &  &  &  &  &  \\ \hline
& $0.01$ &  & $0.121$ &  &  &  &  &  &  \\ \hline
& $0.02$ &  & $0.484$ &  &  &  &  &  &  \\ \hline
& $0.04$ &  & $1.94$ &  &  &  &  &  &  \\ \hline
\end{tabular}%
\end{equation*}

Last, we will obtain a numerical solution of the reaction-diffusion equation
for $a=100,b=1,d=0.001$ which is a reaction dominated case with stiff
reaction.

In Table 5, $L_{2}$ and $L_{\infty }$ error norms are calculated both for $U$
and $V,$ for $N=512$ and various $\Delta t.$

\begin{equation*}
\begin{tabular}{|l|l|l|l|l|l|l|l|l|l|}
\hline
\multicolumn{10}{|l|}{Table 5: Error norms $L_{2}$ and $L_{\infty }$ for
diffusion dominated case with stiff reaction} \\ \hline
\multicolumn{10}{|l|}{for $a=100,$ $b=1,d=0.001$} \\ \hline
\multicolumn{6}{|l}{TQB} & \multicolumn{4}{|l|}{Polynomial quintic B-spline,%
\cite{b15}} \\ \hline
$N$ & $\Delta t$ & $U$ &  & $V$ &  & $U$ &  & $V$ &  \\ \hline
&  & $L_{2}\times 10^{5}$ & $L_{\infty }\times 10^{5}$ & $L_{2}\times 10^{3}$
& $L_{\infty }\times 10^{3}$ & $L_{2}\times 10^{5}$ & $L_{\infty }\times
10^{5}$ & $L_{2}\times 10^{3}$ & $L_{\infty }\times 10^{3}$ \\ \hline
$512$ & $0.005$ & \multicolumn{1}{|r|}{$0,068087$} & \multicolumn{1}{|r|}{$%
0,076753$} & \multicolumn{1}{|r|}{$0,067406$} & \multicolumn{1}{|r|}{$%
0,075986$} & \multicolumn{1}{|r|}{$0,068124$} & \multicolumn{1}{|r|}{$%
0,076795$} & \multicolumn{1}{|r|}{$0,067443$} & \multicolumn{1}{|r|}{$%
0,076027$} \\ \hline
& $0.01$ & \multicolumn{1}{|r|}{$0,272462$} & \multicolumn{1}{|r|}{$0,307141$%
} & \multicolumn{1}{|r|}{$0,269738$} & \multicolumn{1}{|r|}{$0,30407$} & 
\multicolumn{1}{|r|}{$0,272499$} & \multicolumn{1}{|r|}{$0,307183$} & 
\multicolumn{1}{|r|}{$0,269774$} & \multicolumn{1}{|r|}{$0,304111$} \\ \hline
& $0.02$ & \multicolumn{1}{|r|}{$1,089996$} & \multicolumn{1}{|r|}{$1,228729$%
} & \multicolumn{1}{|r|}{$1,079096$} & \multicolumn{1}{|r|}{$1,216442$} & 
\multicolumn{1}{|r|}{$1,090033$} & \multicolumn{1}{|r|}{$1,228771$} & 
\multicolumn{1}{|r|}{$1,079133$} & \multicolumn{1}{|r|}{$1,216484$} \\ \hline
& $0.04$ & \multicolumn{1}{|r|}{$4,360663$} & \multicolumn{1}{|r|}{$4,915684$%
} & \multicolumn{1}{|r|}{$4,317057$} & \multicolumn{1}{|r|}{$4,866527$} & 
\multicolumn{1}{|r|}{$4,360700$} & \multicolumn{1}{|r|}{$4,915725$} & 
\multicolumn{1}{|r|}{$4,317093$} & \multicolumn{1}{|r|}{$4,866568$} \\ \hline
\multicolumn{10}{|l|}{CN-MG method \cite{b3}} \\ \hline
$512$ & $0.005$ &  &  &  & $0.0760$ &  &  &  &  \\ \hline
& $0.01$ &  &  &  & $0.304$ &  &  &  &  \\ \hline
& $0.02$ &  &  &  & $1.22$ &  &  &  &  \\ \hline
& $0.04$ &  &  &  & $4.87$ &  &  &  &  \\ \hline
\end{tabular}%
\end{equation*}

\subsubsection{Nonlinear Problem (Brusselator Model)}

The Brusselator model is a general nonlinear reaction-diffusion system that
models predicting oscillations in chemical reactions. The system was firs
presented by Prigogine and Lefever \cite{b6} showing two variable
autocatalytic reactions. This is one of the simplest reaction-diffusion
equations exhibiting Turing instability, and that large-scale studies have
been conducted on this model with the system being investigated both
analytically and numerically. The general reaction-diffusion equation system
for this model given as:

\begin{equation}
\begin{tabular}{l}
$\dfrac{\partial U}{\partial t}=\varepsilon _{1}\dfrac{\partial ^{2}U}{%
\partial x^{2}}+A+U^{2}V-(B+1)U$ \\ 
$\dfrac{\partial V}{\partial t}=\varepsilon _{2}\dfrac{\partial ^{2}V}{%
\partial x^{2}}+BU-U^{2}V$%
\end{tabular}
\label{r9}
\end{equation}%
where $\varepsilon _{i},i=1,2$ are diffusion constants, $x$ is the spatial
coordinate and $U,V$ are functions of $x$ and $t$ representing
concentrations The initial conditions are selected similar to the reference 
\cite{b9}.

\begin{equation}
\begin{tabular}{ll}
$U(x,0)=0.5,$ & $V(x,0)=1+5x$%
\end{tabular}
\label{r28}
\end{equation}%
and the additional boundary conditions 
\begin{equation*}
\begin{tabular}{ll}
$U_{xx}(x_{0},t)=0$ & $U_{xx}(x_{N},t)=0,$ \\ 
$V_{xx}(x_{0},t)=0$ & $V_{xx}(x_{N},t)=0.$%
\end{tabular}%
\end{equation*}%
In the equation system (\ref{r9}), the coefficients are taken as $\
\varepsilon _{1}=\varepsilon _{2}=10^{-4},$ $A=1,$ $B=3.4.$The solutions are
obtained in the region $x\in \left[ 0,1\right] $, and the programme is run
by the time $t=15$; for space discretization $N=200$ split point and for
time discretization \ $\Delta t=0.01$ time step is used. The solutions under
these selections, are given in Fig. 1 and Fig. 2. which show changes of
density of the functions.\ When wave action is examined, we observe that
both $U$ and $V$ exhibit periodic wave motion under these conditions.

\begin{equation*}
\begin{array}{c}
\FRAME{itbpF}{3.0242in}{2.3999in}{0in}{}{}{figure1.png}{\special{language
"Scientific Word";type "GRAPHIC";maintain-aspect-ratio TRUE;display
"USEDEF";valid_file "F";width 3.0242in;height 2.3999in;depth
0in;original-width 5.073in;original-height 4.0205in;cropleft "0";croptop
"1";cropright "1";cropbottom "0";filename '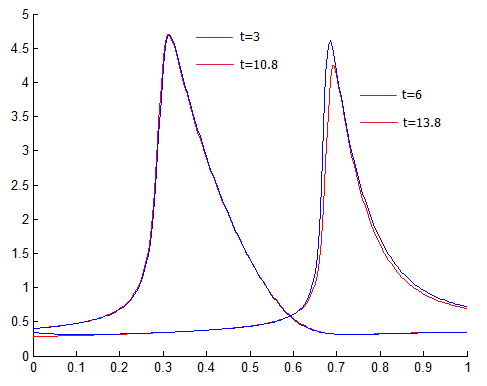';file-properties
"XNPEU";}} \\ 
\text{Figure 1: Periodic wave motion for U} \\ 
\text{for }N=200\text{ }\Delta t=0.01%
\end{array}%
\begin{array}{c}
\FRAME{itbpF}{3.0441in}{2.4042in}{0.0623in}{}{}{figure2.png}{\special%
{language "Scientific Word";type "GRAPHIC";maintain-aspect-ratio
TRUE;display "USEDEF";valid_file "F";width 3.0441in;height 2.4042in;depth
0.0623in;original-width 5.0315in;original-height 3.9686in;cropleft
"0";croptop "1";cropright "1";cropbottom "0";filename
'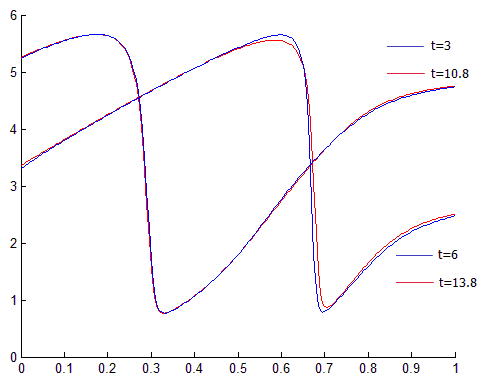';file-properties "XNPEU";}} \\ 
\text{Figure 2: Periodic wave motion for V} \\ 
\text{for }N=200\text{ }\Delta t=0.01%
\end{array}%
\end{equation*}%
The density values for periodic motion are given in Table 6. We see that
this wave is observed as a period of about $7.8$; whereas the period $7.7$
is found when the polynomial quintic B-spline collocation algorithm is
implemented (\ref{r15})

\begin{equation*}
\begin{tabular}{|llllllll|}
\hline
\multicolumn{8}{|l|}{Table 6: Density values for periodic motion for TQB} \\ 
\hline
$Density$ & \multicolumn{1}{|l}{$t$} & \multicolumn{1}{|l}{$x=0.0$} & 
\multicolumn{1}{|l}{$x=0.2$} & \multicolumn{1}{|l}{$x=0.4$} & 
\multicolumn{1}{|l}{$x=0.6$} & \multicolumn{1}{|l}{$x=0.8$} & 
\multicolumn{1}{|l|}{$x=1.0$} \\ \hline
$U$ & \multicolumn{1}{|l}{$3$} & \multicolumn{1}{|r}{0,284595} & 
\multicolumn{1}{|r}{0,317799} & \multicolumn{1}{|r}{0,377380} & 
\multicolumn{1}{|r}{0,604709} & \multicolumn{1}{|r}{1,623703} & 
\multicolumn{1}{|r|}{0,691906} \\ \hline
& \multicolumn{1}{|l}{$10.8$} & \multicolumn{1}{|r}{0,344555} & 
\multicolumn{1}{|r}{0,321243} & \multicolumn{1}{|r}{0,376194} & 
\multicolumn{1}{|r}{0,605486} & \multicolumn{1}{|r}{1,715194} & 
\multicolumn{1}{|r|}{0,716792} \\ \hline
& \multicolumn{1}{|l}{$6$} & \multicolumn{1}{|r}{0,400865} & 
\multicolumn{1}{|r}{0,687572} & \multicolumn{1}{|r}{2,884364} & 
\multicolumn{1}{|r}{0,549937} & \multicolumn{1}{|r}{0,323697} & 
\multicolumn{1}{|r|}{0,348838} \\ \hline
& \multicolumn{1}{|l}{$13.8$} & \multicolumn{1}{|r}{0,398971} & 
\multicolumn{1}{|r}{0,680057} & \multicolumn{1}{|r}{2,911740} & 
\multicolumn{1}{|r}{0,533798} & \multicolumn{1}{|r}{0,322405} & 
\multicolumn{1}{|r|}{0,347582} \\ \hline
&  & \multicolumn{1}{r}{} & \multicolumn{1}{r}{} & \multicolumn{1}{r}{} & 
\multicolumn{1}{r}{} & \multicolumn{1}{r}{} & \multicolumn{1}{r|}{} \\ \hline
$V$ & \multicolumn{1}{|l}{$3$} & \multicolumn{1}{|r}{3,363723} & 
\multicolumn{1}{|r}{4,250910} & \multicolumn{1}{|r}{5,066610} & 
\multicolumn{1}{|r}{5,546754} & \multicolumn{1}{|r}{1,650507} & 
\multicolumn{1}{|r|}{2,507119} \\ \hline
& \multicolumn{1}{|l}{$10.8$} & \multicolumn{1}{|r}{3,309473} & 
\multicolumn{1}{|r}{4,240150} & \multicolumn{1}{|r}{5,062313} & 
\multicolumn{1}{|r}{5,651837} & \multicolumn{1}{|r}{1,591938} & 
\multicolumn{1}{|r|}{2,473710} \\ \hline
& \multicolumn{1}{|l}{$6$} & \multicolumn{1}{|r}{5,258678} & 
\multicolumn{1}{|r}{5,632343} & \multicolumn{1}{|r}{1,073700} & 
\multicolumn{1}{|r}{2,739517} & \multicolumn{1}{|r}{4,300681} & 
\multicolumn{1}{|r|}{4,755329} \\ \hline
& \multicolumn{1}{|l}{$13.8$} & \multicolumn{1}{|r}{5,241915} & 
\multicolumn{1}{|r}{5,634312} & \multicolumn{1}{|r}{1,065232} & 
\multicolumn{1}{|r}{2,769906} & \multicolumn{1}{|r}{4,269058} & 
\multicolumn{1}{|r|}{4,737755} \\ \hline
\end{tabular}%
\end{equation*}

\begin{equation*}
\begin{tabular}{|llllllll|}
\hline
\multicolumn{8}{|l|}{Table 7: Density values for periodic motion for quintic
B-spline \cite{b15}} \\ \hline
Density & \multicolumn{1}{|l}{$t$} & \multicolumn{1}{|l}{$x=0.0$} & 
\multicolumn{1}{|l}{$x=0.2$} & \multicolumn{1}{|l}{$x=0.4$} & 
\multicolumn{1}{|l}{$x=0.6$} & \multicolumn{1}{|l}{$x=0.8$} & 
\multicolumn{1}{|l|}{$x=1.0$} \\ \hline
$U$ & \multicolumn{1}{|l}{$3$} & \multicolumn{1}{|l}{0,284657} & 
\multicolumn{1}{|l}{0,317966} & \multicolumn{1}{|l}{0,377959} & 
\multicolumn{1}{|l}{0,612881} & \multicolumn{1}{|l}{1,519483} & 
\multicolumn{1}{|l|}{0,648434} \\ \hline
& \multicolumn{1}{|l}{$10.7$} & \multicolumn{1}{|l}{0,347747} & 
\multicolumn{1}{|l}{0,321168} & \multicolumn{1}{|l}{0,376204} & 
\multicolumn{1}{|l}{0,611218} & \multicolumn{1}{|l}{1,626310} & 
\multicolumn{1}{|l|}{0,680742} \\ \hline
& \multicolumn{1}{|l}{$6$} & \multicolumn{1}{|l}{0,401741} & 
\multicolumn{1}{|l}{0,706734} & \multicolumn{1}{|l}{2,716642} & 
\multicolumn{1}{|l}{0,510302} & \multicolumn{1}{|l}{0,326204} & 
\multicolumn{1}{|l|}{0,352411} \\ \hline
& \multicolumn{1}{|l}{$13.7$} & \multicolumn{1}{|l}{0,398904} & 
\multicolumn{1}{|l}{0,691408} & \multicolumn{1}{|l}{2,769059} & 
\multicolumn{1}{|l}{0,500480} & \multicolumn{1}{|l}{0,324523} & 
\multicolumn{1}{|l|}{0,350579} \\ \hline
&  &  &  &  &  &  &  \\ \hline
$V$ & \multicolumn{1}{|l}{$3$} & \multicolumn{1}{|l}{3,363896} & 
\multicolumn{1}{|l}{4,251219} & \multicolumn{1}{|l}{5,066734} & 
\multicolumn{1}{|l}{5,537413} & \multicolumn{1}{|l}{1,732740} & 
\multicolumn{1}{|l|}{2,580615} \\ \hline
& \multicolumn{1}{|l}{$10.7$} & \multicolumn{1}{|l}{3,299664} & 
\multicolumn{1}{|l}{4,233913} & \multicolumn{1}{|l}{5,056668} & 
\multicolumn{1}{|l}{5,637796} & \multicolumn{1}{|l}{1,659946} & 
\multicolumn{1}{|l|}{2,534846} \\ \hline
& \multicolumn{1}{|l}{$6$} & \multicolumn{1}{|l}{5,257254} & 
\multicolumn{1}{|l}{5,606791} & \multicolumn{1}{|l}{1,137215} & 
\multicolumn{1}{|l}{2,825295} & \multicolumn{1}{|l}{4,355469} & 
\multicolumn{1}{|l|}{4,798749} \\ \hline
& \multicolumn{1}{|l}{$13.7$} & \multicolumn{1}{|l}{5,234725} & 
\multicolumn{1}{|l}{5,613815} & \multicolumn{1}{|l}{1,119445} & 
\multicolumn{1}{|l}{2,846165} & \multicolumn{1}{|l}{4,317357} & 
\multicolumn{1}{|l|}{4,774541} \\ \hline
\end{tabular}%
\end{equation*}

\subsubsection{\textbf{Nonlinear Problem (Schnakenberg Model)}}

The Schnakenberg model is a well-known reaction-diffusion model which is a
simplified version of the Brusselator model. It is a relatively easy system
for modeling the reaction-diffusion mechanism. There are many studies in the
literature on this model. Firstly it is modeled by Schakenberg \cite{b4} and
given as:

\begin{equation}
\begin{tabular}{l}
$\dfrac{\partial U}{\partial t}=\dfrac{\partial ^{2}U}{\partial x^{2}}%
+\gamma (a-U+U^{2}V)$ \\ 
$\dfrac{\partial V}{\partial t}=d\dfrac{\partial ^{2}V}{\partial x^{2}}%
+\gamma (b-U^{2}V)$%
\end{tabular}
\label{r10}
\end{equation}%
where U and V denote the concentration of activator and inhibitor
respectively, $d$ is diffusion coefficient, $\gamma $, $a$ and $b$ are rate
constants of the biochemical reactions. The oscillation problem is taken
into account for the Schnakenberg Model. Accordingly, the parameters for
system (\ref{r10}) are selected as $a=0.126779$,$b=0.792366,d=10$ and $%
\gamma =10^{4}.$The problem's initial conditions:

\begin{eqnarray}
U(x,0) &=&0.919145+0.001\underset{j=1}{\overset{25}{\sum }}\frac{\cos (2\pi
jx)}{j}  \label{r30} \\
V(x,0) &=&0.937903+0.001\underset{j=1}{\overset{25}{\sum }}\frac{\cos (2\pi
jx)}{j}  \notag
\end{eqnarray}%
are on the interval $[-1.1]$. The boundary conditions left, right and
additional boundary conditions are:%
\begin{equation*}
\begin{tabular}{ll}
$U_{x}(x_{0},t)=0$ & $U_{x}(x_{N},t)=0,$ \\ 
$V_{x}(x_{0},t)=0$ & $V_{x}(x_{N},t)=0.$%
\end{tabular}%
\end{equation*}%
\begin{equation*}
\begin{tabular}{ll}
$U_{xxx}(x_{0},t)=0$ & $U_{xxx}(x_{N},t)=0,$ \\ 
$V_{xxx}(x_{0},t)=0$ & $V_{xxx}(x_{N},t)=0.$%
\end{tabular}%
\end{equation*}%
Computations are performed to the $t=2.5$ for space/time combinations given
in Table. 8. Obtained relative error values are given in Table 8 together
with the results of quintic B-spline collocation method \cite{b15}.

\begin{equation*}
\begin{tabular}{|c|c|c|c|c|c|}
\hline
\multicolumn{6}{|c|}{Tablo 8: Relative error values for $N=100$ in $t=2.5$}
\\ \hline
$\Delta t$ & Nu. of steps & $U$ & $U$\cite{b15} & $V$ & $V$\cite{b15} \\ 
\hline
$5\times 10^{-6}$ & 500000 & $0$ & $5.7160\times 10^{-14}$ & $5.4418\times
10^{-17}$ & $5.4564\times 10^{-14}$ \\ \hline
$5\times 10^{-5}$ & 50000 & $6.2202\times 10^{-17}$ & $1.5653\times 10^{-10}$
& $1.6794\times 10^{-16}$ & $1.1105\times 10^{-10}$ \\ \hline
$1\times 10^{-4}$ & 25000 & $1.7593\times 10^{-16}$ & $9.8744\times 10^{-10}$
& $2.4423\times 10^{-16}$ & $8.8599\times 10^{-10}$ \\ \hline
$1.20\times 10^{-4}$ & 20833 & $1.5668\times 10^{-16}$ & $1.5055\times
10^{-09}$ & $2.2996\times 10^{-16}$ & $1.3790\times 10^{-09}$ \\ \hline
$1.32\times 10^{-4}$ & 18939 & $1.4610\times 10^{-16}$ & $1.0564\times
10^{-01}$ & $2.9664\times 10^{-16}$ & $1.0301\times 10^{-01}$ \\ \hline
$1\times 10^{-3}$ & 2500 & $2.5895\times 10^{-14}$ & - & $2.0341\times
10^{-14}$ & - \\ \hline
$2\times 10^{-3}$ & 1250 & $5.4591\times 10^{-09}$ & - & $3.9448\times
10^{-09}$ & - \\ \hline
$5\times 10^{-3}$ & 500 & $5.4960\times 10^{-06}$ & - & $4.7003\times
10^{-06}$ & - \\ \hline
\end{tabular}%
\end{equation*}

As can be seen from Table 8, the algorithm produces accurate results even
when the time increment is larger\ The Figure 3 was drawn to show the
oscillation movements for values $\Delta t=5\times 10^{-5}$, $N=100$ and $%
N=200$ It is shown in Fig. 3 that the functions $U$ and $V$ make 9
oscillations when $N=200$ and $N=100.$This result with the references \cite%
{b10} and \cite{b11} shows that a finer mesh is necessary for accurate
solutions.

\begin{equation*}
\begin{tabular}{l}
\FRAME{itbpF}{6.3771in}{2.8297in}{0in}{}{}{figure3.png}{\special{language
"Scientific Word";type "GRAPHIC";maintain-aspect-ratio TRUE;display
"USEDEF";valid_file "F";width 6.3771in;height 2.8297in;depth
0in;original-width 10.6458in;original-height 4.708in;cropleft "0";croptop
"1";cropright "1";cropbottom "0";filename '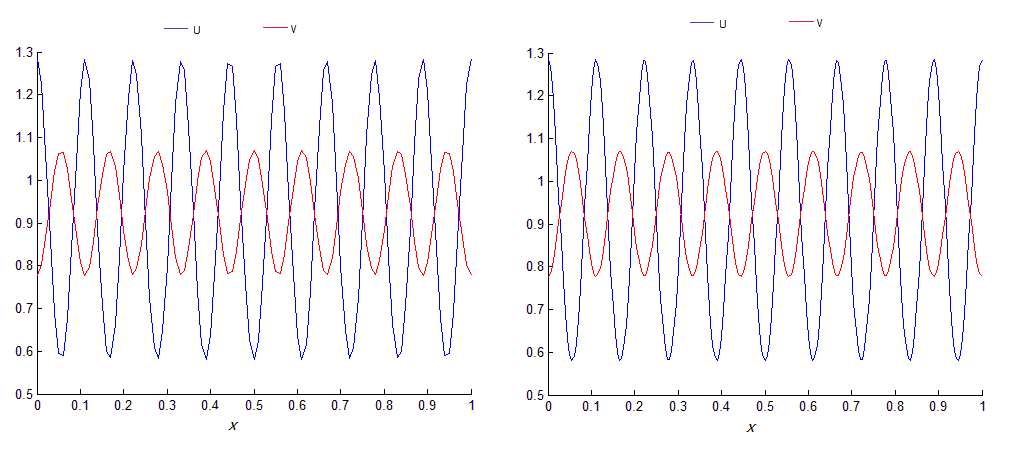';file-properties
"XNPEU";}} \\ 
Fig. 3: The oscillation movement for $N=100$ and $N=200$ in the moment $%
t=2.5 $%
\end{tabular}%
\end{equation*}

\subsubsection{\textbf{Nonlinear Problem (Gray-Scott Model)}}

The Gray-Scott model is a reaction-difussion system which models the forming
of certain spatial patterns by a few chemical species, that exit in the
nature. It was put forward by Gray and Scott \cite{b2} and the
reaction-diffusion system is given as:

\begin{equation}
\begin{tabular}{l}
$\dfrac{\partial U}{\partial t}=\varepsilon _{1}\dfrac{\partial ^{2}U}{%
\partial x^{2}}-U^{2}V+f(1-U)$ \\ 
$\dfrac{\partial V}{\partial t}=\varepsilon _{2}\dfrac{\partial ^{2}V}{%
\partial x^{2}}+U^{2}V-(f+k)V$%
\end{tabular}
\label{r11}
\end{equation}%
In this section, the numerical method was tested with repetitive spot
patterns on Gray-Scott Model. The parameters for system (\ref{r11}) were
chosen as the reference \cite{b22}%
\begin{equation*}
\begin{tabular}{llll}
$\varepsilon _{1}=1,$ & $\varepsilon _{2}=0.01,$ & $a=9$ & $b=0.4$%
\end{tabular}%
\end{equation*}%
with these parameters the initial conditions of system (\ref{r11}) were
taken as

\begin{equation}
\begin{tabular}{l}
$U(x,0)=1-\frac{1}{2}\sin ^{100}(\pi \frac{(x-L)}{2L})$ \\ 
$V(x,0)=\frac{1}{4}\sin ^{100}(\pi \frac{(x-L)}{2L})$%
\end{tabular}
\label{r31}
\end{equation}%
and solutions were investigated in interval $[-L,L]$ \ and $L=50$. For space
discretization $N=400$ and for time discretization $\Delta t=0.2$ were
selected. Dirichlet boundary condititions

\begin{equation*}
\begin{tabular}{l}
$U(x_{0},t)=U(x_{N},t)=1,$ \\ 
$V(x_{0},t)=V(x_{N},t)=0$%
\end{tabular}%
\end{equation*}%
together with additonal Neuman boundary conditions

\begin{equation*}
\begin{tabular}{l}
$U_{x}(x_{0},t)=U_{x}(x_{N},t)=0,$ \\ 
$V_{x}(x_{0},t)=V_{x}(x_{N},t)=0$%
\end{tabular}%
\end{equation*}%
are used. Numerical computations were made until $t=100$ \ and $t=500$ so
that repetitive patterns were obtained. Under these initial conditions,
primarily two pulses were created \ and separated from each other, with each
pulse then being split into two again to form four pulses, as shown in Fig.
5. until time $t=1000,$ as time evolved. This self-replicating process goes
on to cover the spatial domain. These splitting movements of the functions $%
U $ and $V$ due to time and space are presented in Figs 4-5.

\begin{equation*}
\begin{tabular}{l}
\begin{tabular}{l}
\FRAME{itbpF}{6.7205in}{2.802in}{0in}{}{}{figure4.png}{\special{language
"Scientific Word";type "GRAPHIC";maintain-aspect-ratio TRUE;display
"USEDEF";valid_file "F";width 6.7205in;height 2.802in;depth
0in;original-width 10.5308in;original-height 4.3751in;cropleft "0";croptop
"1";cropright "1";cropbottom "0";filename 'Figure4.png';file-properties
"XNPEU";}}%
\end{tabular}
\\ 
Figure 4: The splitting process of repetitive spot pattern of waves for $%
t=100$ and $t=500$%
\end{tabular}%
\end{equation*}%
\begin{equation*}
\begin{tabular}{l}
\begin{tabular}{l}
\  \  \  \  \  \  \  \  \  \  \  \  \  \  \  \  \  \  \  \  \FRAME{itbpF}{3.1946in}{2.4621in}{0in%
}{}{}{figure5.png}{\special{language "Scientific Word";type
"GRAPHIC";display "USEDEF";valid_file "F";width 3.1946in;height
2.4621in;depth 0in;original-width 5.8332in;original-height 4.3751in;cropleft
"0";croptop "1";cropright "1";cropbottom "0";filename
'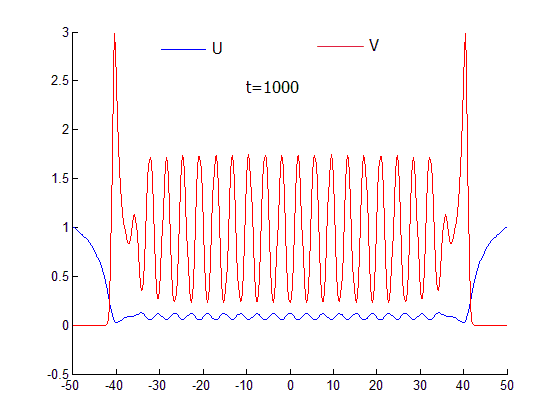';file-properties "XNPEU";}}%
\end{tabular}
\\ 
Figure 5: The splitting process of repetitive spot pattern of waves for $%
t=1000$%
\end{tabular}%
\end{equation*}%
The intensity changes of functions $U$ and $V$ due to time and space are
presented in Fig. 6 and Fig. 7, respectively. These spatial patterns, which
known as repetitive spot patterns, initially starting with two waves of
splitting movement, seem to cover the whole domain with branching over time.

\begin{equation*}
\begin{tabular}{l}
\FRAME{itbpF}{3.5587in}{3.2344in}{0in}{}{}{figure6.png}{\special{language
"Scientific Word";type "GRAPHIC";display "USEDEF";valid_file "F";width
3.5587in;height 3.2344in;depth 0in;original-width 7.0136in;original-height
4.875in;cropleft "0";croptop "1";cropright "1";cropbottom "0";filename
'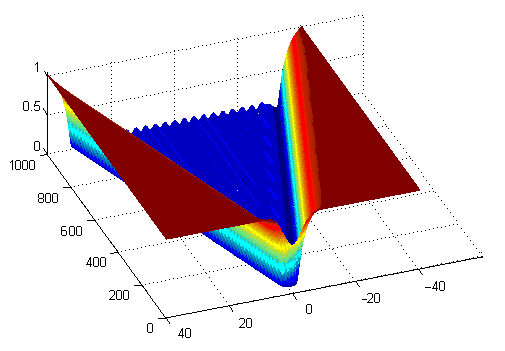';file-properties "XNPEU";}} \\ 
Fig. 6: Repetitive spot pattern of waves of the function U%
\end{tabular}%
\end{equation*}

\begin{equation*}
\begin{tabular}{l}
\FRAME{itbpF}{3.6729in}{3.1021in}{0in}{}{}{figure7.png}{\special{language
"Scientific Word";type "GRAPHIC";display "USEDEF";valid_file "F";width
3.6729in;height 3.1021in;depth 0in;original-width 7.542in;original-height
5.4864in;cropleft "0";croptop "1";cropright "1";cropbottom "0";filename
'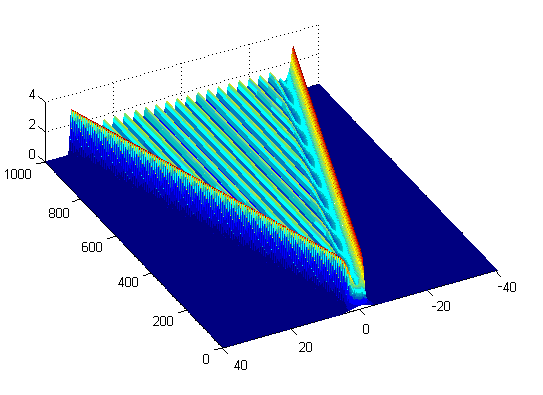';file-properties "XNPEU";}} \\ 
Fig. 7: Repetitive spot pattern of waves of the function V%
\end{tabular}%
\end{equation*}

\subsection{Discussion}

Proposed algorithm has been used for calculating numerical solutions of
reaction-diffusion equation systems. Solutions of linear and nonlinear RD
systems are shown on the models of certain chemical problems: the
Brusselator model, Schnakenberg model and Gray-Scott models are simulated
suitably by use of the suggested algorithm. The proposed TQB algorithm is an
alternative method to the more usual polynomial quintic B-spline collocation
methods (PQBCM). The results of the suggested algorithm are documented
together with those obtained with PQBCM and Crank-Nicolson Multigrid solver
method (CN-MG) for the test problem. It can be seen from the tables (3-5)
that the accuracy of the algorithms are almost the same with that for the
PQBCM and are better than the CN-MG. Solutions of the nonlinear problems,
which have no analytical solutions in general, are given graphically. Model
solutions are represented fairly and can be compared with the equivalent
graphs given in the studies \cite{b9,b10,b11,b15}. Use of the trigonometric
B-spline having continuity of the order four allow us to have an approximate
functions in the order of four. Therefore, differential equations in the
order of four can be solved numerical by using the trigonometric B-spline
functions to have solutions of continuity in the order of four.
Consequently, the TQB collocation method produces fairly acceptable results
for numerical solutions of reaction-diffusion systems. Thus, it is also
recommended to finding solutions of the other partial dfferential equations.

\textbf{Competing Interest}

The authors declare that they have no competing interests.

\textbf{Authors' contributions}

AOT carried out the algorithm implementation and conducted test studies,
participated in the sequence alignment and drafted the manuscript. NA
conceived of the study, participated in developing algorithm for Numeric
Solutions of trigonometric quintic B-splines and helped to draft the
manuscript. ID participated in its design and coordination and helped to
draft the manuscript. All authors read and approved the final manuscript.

\textbf{Acknowledgement}

The paper is presented at the conference in Ankara, International Conference
on Applied Mathematics and Analysis (ICAMA2016).


\begin{thebibliography}{99}
\bibitem{b1} Nikolis A., Numerical solutions of ordinary differential
equations with quadratic trigonometric splines,\ Applied Mathematics
E-Notes, 4(1995),142-149.

\bibitem{b2} Gray P. and Scott S.K., Autocatalytic reactions in the
isothermal, continuous stirred tank reactor: oscillations and instabilities
in the system A+2B 3B, B C, Chem. Eng.Sci., 39(1984), 1087-1097.

\bibitem{b3} Chou C., Zhang Y.,Zhao R.and Nie Q., Numerical methods for
stiff reaction-diffusion systems, Discrete and Continuous Dynamical
Systems-Series B, 7(2007), 515-525.

\bibitem{b4} Schnakenberg, J.Simple chemical reaction systems with limit
cycle behavior, J.Theoret. Biol., 81(1979), 389-400.

\bibitem{b6} Prigogine, I. and Lefever, R. Symmetry breaking instabilities
in dissipative systems, J. Chem. Phys., 48(1968), 1695-1700.

\bibitem{b8} Nikolis, A. and Seimenis, I. Solving dynamical systems with
cubic trigonometric splines, Applied Mathematics E-notes, 5(2005), 116-123.

\bibitem{b9} Zegeling, P.A and Kok, H.P. Adaptive moving mesh computations
for reaction-diffusion systems, Journal of Computational and Applied
Mathematics, 168(2004), 519-528.

\bibitem{b10} Madzvamuse, A. Wathen, A.J. and Maini, P.K. A moving grid
finite element method applied to a biological pattern generator, Journal of
Computational Physics, 190(2003), 478-500.

\bibitem{b11} Ruuth, S.J. Implicit-explicit methods for reaction-diffusion
problems in pattern formation, Journal of Mathematical Biology, 34(1995),
148-176.

\bibitem{b15} Sahin, A: Numerical solutions of the reaction-diffusion
equations with B-spline finite element method. Dissertation, Eski\c{s}ehir
Osmangazi University, Eski\c{s}ehir, Turkey (2009).

\bibitem{b16} Hamid, N.N.A. Majid, A.A. and Ismail, A.I.M. Cubic
trigonometric B-Spline applied to linear two-point boundary value problems
of order,World Academy of Science, Engineering and Technology, 70(2010),
798-803.

\bibitem{b18} Gupta, Y.and Kumar, M. A computer based numerical method for
singular boundary value problems, International Journal of Computer
Applications,30(1)(2011), 21-25,.

\bibitem{b19} Abbas, M. Majid, A. A. Ismail, A. I. M and Rashid, A. The
application of the cubic trigonometric B-spline to the numerical solution of
the hyperbolic problems, Applied Mathematica and Computation, 239(2014)
74-88.

\bibitem{b20} Abbas, M. Majid, A. A. Ismail, A. I. M and Rashid, A.
Numerical method using cubic trigonometric B-spline tecnique for
nonclassical diffusion problems, Abstract and applied analysis,
2014(2014),1-12.

\bibitem{kgo} Zin, S. M. Abbas, M. Majid, A.A and Ismail, A.I. M. A new
trigonometric spline approach to numerical solution of generalized nonlinear
Klien-Gordon equation, PLOS one, 9(5)(2014),1-9

\bibitem{b21} Schoenberg, I. J. On trigonometric spline interpolation, J.
Math. Mech.13(5)(1964), 795--825.

\bibitem{b22} Craster, R.V. and Sassi, R. Spectral algorithms for
reaction-diffusion equations,Technical Report. Note del Polo, No. 99 (2006).
\end{thebibliography}
\end{document}